# Integrated PV Charging of EV Fleet Based on Dynamic Energy Prices and Offer of Reserves

Gautham Ram Chandra Mouli, *Student Member, IEEE*, Mahdi Kefayati, *Member, IEEE*,
Ross Baldick, *Fellow, IEEE*, Pavol Bauer, *Senior Member, IEEE*

*Abstract*—Workplace charging of electric vehicles (EV) from photovoltaic (PV) panels installed on an office building can provide several benefits. This includes the local production and use of PV energy for charging the EV and making use of dynamic tariffs from the grid to schedule the energy exchange with the grid. The long parking time of EV at the workplace provide the chance for the EV to support the grid via vehicle-to-grid technology, the use of a single EV charger for charging several EV by multiplexing and the offer of ancillary services to the grid for up and down regulation. Further, distribution network constraints can be considered to limit the power and prevent the overloading of the grid. A single MILP formulation that considers all the above applications has been proposed in this paper for a charging a fleet of EVs from PV. The MILP is implemented as a receding-horizon model predictive energy management system. Numerical simulation based on market and PV data in Austin, Texas have shown 31% to 650% reduction in the cost of EV charging when compared to immediate and average rate charging policies.

## I. Introduction

Electric vehicles (EV) provide a zero-emission, low noise and highly efficient mode of transportation. The current estimate for the USA is that there will be 1.2 million EV by 2020 [1]. Electric vehicles are, however, sustainable only if the electricity used to charge them comes from sustainable sources. Electricity generated from a fuel mix that is largely dominated by fossil fuels does not eliminate the emissions but mostly moves it from the vehicle to the power plant [2], [3]. While this can have environmental advantages, complete elimination of emissions is contingent on utilizing non-emitting resources for electricity production. It is here that the phenomenal growth in the use of photovoltaic (PV) system for distributed generation and its falling cost over the years can have a direct impact.

EVs used to commute to work are parked at the workplace for long hours during the day when the sun is shining. Workplaces like industrial sites and office buildings harbor a great potential for PV panels with their large surfaces on flat

Gautham Ram and Pavol Bauer are with the Department of Electrical Sustainable Energy, Delft University of Technology, 2628 CD Delft, The Netherlands. (e-mail: P.Bauer@tudelft.nl, G.R.Chandamouli@tudelft.nl).

Ross Baldick and Mahdi Kefayati are with the Department of Electrical and Computer Engineering, The University of Texas at Austin, Austin, TX 78712. (e-mail: baldick@ece.utexas.edu, kefayati@utexas.edu).

This work was supported by the TKI Switch2SmartGrids grant, Netherlands and the Partners for International Business(PIB) program of the Rijksdienst voor Ondernemend (RVO), Netherlands.

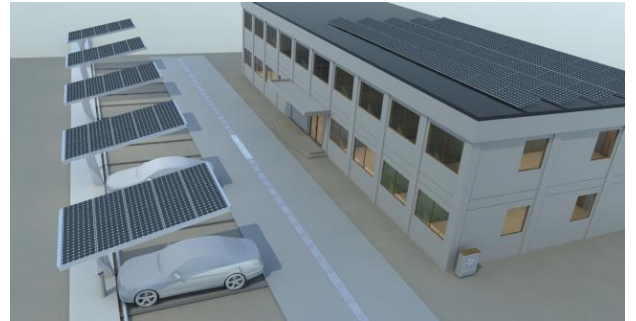

Fig. 1. Design of solar powered EV charging station

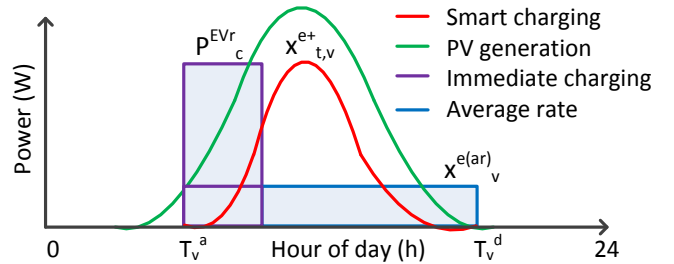

Fig. 2. Immediate, average rate and smart charging of EV

roofs. This potential is largely unexploited today. Energy generated from PV array installed at the workplace can hence be used for charging EVs as shown in Fig. 1. This has several benefits namely:
1. EV battery doubles up as an energy storage for the PV
2. Negative impact of large-scale PV and EV integration on distribution network is mutually reduced [4], [5]
3. Long parking time of EV paves way for implementation of Vehicle-to-grid (V2G) technology where the EV can offer energy and ancillary services to the grid [6], [7].
4. Cost of EV charging from solar is cheaper than charging from the grid and net $CO_2$ emission is lowered [2], [8].

### A. Immediate and average rate charging

Today, when an EV arrives at the workplace and is connected to an electric vehicle supply equipment (EVSE), the EV typically starts charging immediately at the nominal EVSE power rating $P_c^{EVr}$. The charging continues at a constant power till the battery is full[1]. This is referred to as immediate charging (IMM) or uncontrolled charging [9]. This is the simplest form of charging requiring no information from the user or communication infrastructure and results in the lowest

2charging time. However, IMM typically results in a huge demand on the grid, depending on the EVSE capability, as shown in Fig. 2.

EVs parked at the workplace usually have long parking times and this offers the flexibility in scheduling the charging in terms of both charging power and duration. This means that EVs can be charged at a much lower power than the EVSE nominal rating if the EV user arrival time, $T_v^a$, departure time, $T_v^d$ and required energy demand, $d_v$ are known. One approach is the "Average Rate" (AR) charging policy [9].

$$x_v^{e(ar)} = Min.\left\{\frac{d_v}{T_v^d - T_v^a} \le x_v^{ub}, P_c^{EVr}\right\} \quad \forall \, t \in \{T_v^a, T_v^d\} \quad (1)$$

Here, the charging power $x_v^{e(ar)}$ is the minimum of the EVSE capacity and the ratio of the energy demand divided by the parking time of the EV[1]. The advantage of the AR policy is that the charging of the fleet is spread through the day instead of being concentrated in the arrival time (typically early morning), as seen in Fig. 2.

However, both IMM and AR strategies are not 'smart' as it has no correlation to local renewable generation, distribution network capacity constraints and/or energy prices.

*B. Smart charging*

The optimal way to charge EVs is hence to use an energy management system (EMS) that can schedule the charging of an EV fleet by taking into consideration the EV user preferences, local renewable generation and energy prices from the market. Fig. 2 shows the example of a smart charging where the EV charging follows the PV generation. Further, EVs can have extremely fast ramp up and ramp down rates. Chademo and Combo EV charging standards for DC charging stipulate response time of 200ms [10]. This makes EVs ideal candidates for providing ancillary services in the form of reserve capacity to the grid [6], [7], [11], [12].

Following the formulation in [12], [13], an Energy Services Company (ESCo) company acts as an intermediary between the wholesale market operated by the Independent System Operator (ISO) and the EV end-users. The ESCo operates at the workplace where employees drive to the office with an EV and the building has overhead PV installation or a solar carport. The motive of the ESCo is to schedule the charging of the EV and feeding of PV power to the grid in such a way that charging costs are lowered, regulation services are offered to the ISO and at the same time, the income from PV are increased. ESCo achieves this motive by using an EMS to schedule the EV based on a multitude of inputs namely:

1. Information from the EV user about arrival and departure times, the state of charge (SOC) and energy demand.
2. $p_t^{e(buy)}, p_t^{e(sell)}$ are the settlement point prices for buying and selling electricity from the grid at time *t*.
3. $p_t^{r(up)}, p_t^{r(dn)}$ are the clearing prices for capacity for offering reserves to the ISO for up and down regulation.
4. Distribution network constraints $P_t^{DN+}, P_t^{DN-}$ which are the upper limits for drawing and feeding power between the EV car park and the grid at time *t*. These values can be adjusted to implement demand side management (DSM).
5. Solar forecast information: $P_t^{PV(fc)}$ is the generation forecast for a 1kW$_p$ PV system installed at the workplace. The use of solar forecast data in the EMS will help in reducing the uncertainties due to variability in PV generation on diurnal and seasonal basis [14].

*C. Literature review and overview of contributions*

Several earlier works have formulated the optimization problem to charge EV based on renewable generation, energy prices and offer of ancillary services.

Fuzzy logic is used to optimize the EV charging based on PV generation forecast, energy prices in [15] and V2G frequency regulation, grid energy exchange in [16]. The disadvantage is that the use of fuzzy logic without optimization techniques does not guarantee that the obtained solution is optimal.

In [12], [13], linear programming (LP) is used to find the optimal EV strategy for charging and offering reserves based on market prices. In [17], LP is used to reduce the cost of charging EV from PV based on time of use tariffs and PV forecasting. Cost reduction of 6% and 15.2% compared to the base case are obtained for simulation for 12 EV powered from a 50kW PV system. The LP formulation in [18] and heuristic methods used in [19] aim to achieve the two goals: increasing the PV self-consumption in a micro-grid by charging of EVs and reducing the dependency on the grid. However, there is no consideration for time of use tariffs without which there is no incentive to achieve the two goals. In [20], LP is used for planning EV charging based on renewable power forecasting, spinning reserve and EV user requirements in a micro-grid

Stochastic programming is used in [21] to charge EV and offer regulation services based on day-ahead and intraday market prices. For a case study with 50 EV, cost reduction in the tune of 1% to 15% was achieved.

Main contributions of this work include:
- Proposing a single comprehensive model that captures charging of EV from PV, use of dynamic grid prices, implementation of V2G for grid support, using EV to offer ancillary services and considering distribution network capacity constraints as a single MILP problem. Earlier works have considered these as separate optimization problems or as a combination of two or three applications. The disadvantage of the earlier approach is that each optimization problem gives a different optimized EV charging profile and all these profiles cannot be implemented on the same EV at the same time! The best approach is to combine them into one formulation which will then yield a single optimized EV charging profile
- Demonstrating that the formulation of the abovementioned five aspects into one MILP formulation results in large cost savings, which is much higher than what has been achieved earlier.
- With a large number of EV parked at the workplace with long parking times, multiplexing a few EVSE to a larger

---
[1] The expression does not consider the duration in the constant-voltage (CV) charging mode, which occurs typically when EV battery is above 80% SOC and the maximum charging power is limited [28].

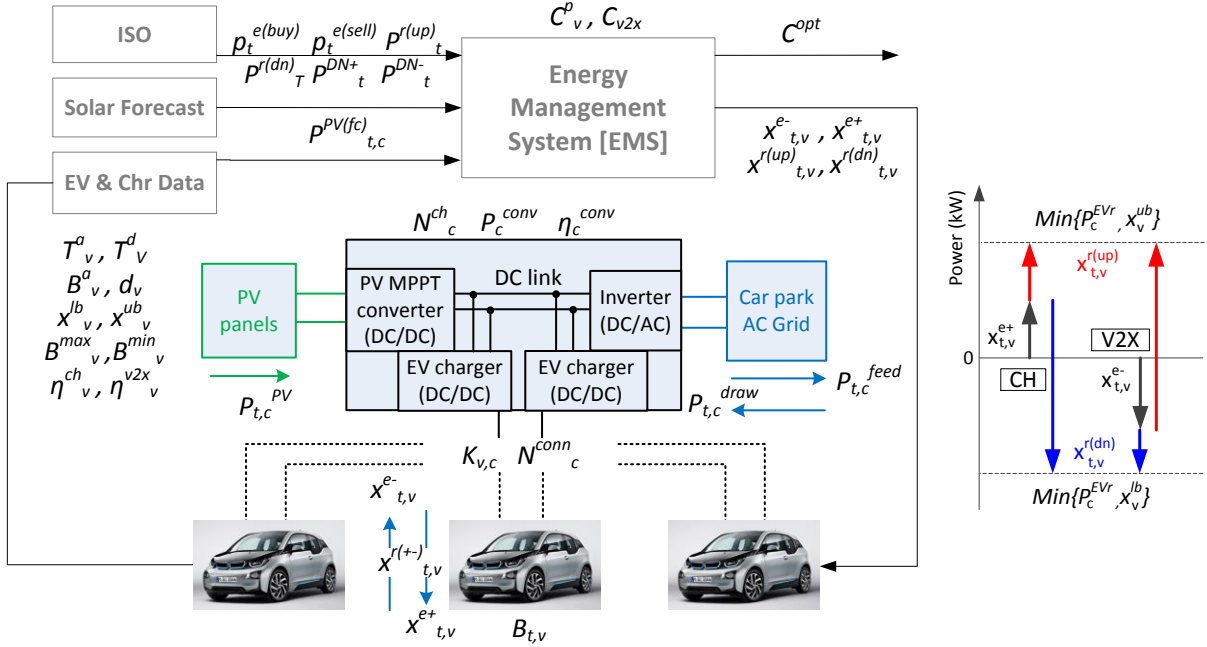

Fig. 3. (Left) Schematic of the Energy Management System (EMS) for the solar powered EV parking garage. $N_c^{ch}$ of the total $N_c^{conn}$ EV connected to each EV-PV charger can be simultaneously charged or discharged, where $N_c^{ch} \leq N_c^{conn}$ (Right) offer of reserve power capacity $x_{t,v}^{r(up)}, x_{t,v}^{r(dn)}$ for up and down regulation during charging (CH) and discharging (V2G) of EV.

number of EVs is a cost effective strategy [22], [23]. These EVSE could offer simultaneous charging of several EV or charge one EV at a time. The scheduling of the multiplexing is formulated for the first time in this work using an MILP formulation.

- Implementation of the optimization using C# code, SQL server and Microsoft solver foundation making it ready for hardware implementation.

*D. Structure of the paper*

Section II describes the layout and parameters of the EMS and the EV-PV car park infrastructure. In section III, the MILP formulation of the EMS is explained and the parameters, constraints and objective function are elaborated. Section IV uses PV generation data and market data for Austin, TX to estimate the optimized cost of charging an EV fleet from PV. The costs are compared to immediate and average rate charging policies to evaluate the cost reduction.

## II. PRELIMINARIES AND INPUTS

*A. Layout of the EMS*

The schematic of the EV-PV charger and the EMS used by the ESCo to optimize the EV charging is shown in Fig. 3.

*1) EV and user input*

If $v$ is the index of the EV and total number of EVs is $V$, each EV arrives at the car park with a state of charge $B_v^a$ at time $T_v^a$ and is parked at one of the several EV-PV chargers. The EV owners provide the information to the EMS about their expected departure time $T_v^d$ and charging energy demand $d_v$. This means that the departure SOC of the vehicle $B_v^d$ is:

$$B_v^d = B_v^a + d_v \quad (2)$$

If the required SOC is not reached by the departure time, the EV owner will be compensated by the ESCo at the rate of $C_v^p$ \$/kWh. The users can enter the maximum and minimum allowed SOC of the EV $B_v^{min}, B_v^{max}$ and the maximum charging and discharging power $x_v^{ub}, x_v^{lb}$ respectively. By setting $x_v^{lb}$ to a non-zero value, the users can choose to participate in V2G services. The efficiency of the EV battery for charging and discharging is $\eta_v^{ch}, \eta_v^{v2x}$ and is either obtained from the EV or stored in a database within the EMS for different EV models.

*2) EV-PV charger*

The 'EV-PV charger' as the term used here is an integrated power converter that consists of three ports to connect to the EVs, PV and the AC grid, as shown in Fig. 3 [14], [22]. Each EV-PV charger is connected to a PV array of rated power $P_c^{PVr}$ via a maximum power point tracking (MPPT) DC/DC converter [24]. The output of the DC/DC PV converter is connected to an internal DC-link. The DC-link is connected to the grid via a DC/AC inverter of rated power $P_c^{conv}$, such that $P_c^{PVr} \leq P_c^{conv}$. There are $N_c^{ch}$ number of isolated DC/DC converter for EV charging that are connected to the DC-link and each have a rated power $P_c^{EVr}$. All power exchanges between any of the three ports namely PV, EV, grid happens via the DC-link.

This integrated converter provides several benefits compared to using separate converters for PV and EV connected over the 50Hz AC grid. First, direct interconnection of the PV and EV over a DC-link is more efficient than an AC interconnection [25], [26]. Second, the integrated converter requires one common inverter to the AC grid instead of separate inverters for PV and EV. This reduces the component count and cost of the converter [22]. Third, by making the isolated DC/DC converter for the EV bidirectional, the EV can now offer V2G services via the integrated converter.

Due to the long parking times of EVs at the workplace, it is economical to use a single EVSE that can be multiplexed to several EVs, with the possibility to charge the EVs simultaneously or sequentially as shown in Fig. 3. Therefore, $N_c^{conn}$ EVs can be connected to each EV-PV charger via DC isolators. The binary variable $K_{v,c} = 1$ indicates the physical connection of $v^{th}$ EV with $c^{th}$ charger and a zero value indicates otherwise.

Each EV-PV charger has $N_c^{ch}$ number of isolated DC/DC converters, where $N_c^{ch} \leq N_c^{conn}$. As per the EV charging standards [27], each EV must be connected to separate power converter and isolated from all power sources. This means that $N_c^{ch}$ of the total $N_c^{conn}$ EVs connected to each EV-PV charger can be simultaneously charged or discharged. In the simple case where $N_c^{ch} = 1$, $N_c^{conn}=2$ and $P_c^{conv} = P_c^{EVr}$, two EVs are connected to one EV-PV charger and one out of the two can (dis)charge at any time up to a power of $P_c^{conv}$. The binary variable $a_{t,v}^c$ indicates which of the $N_c^{conn}$ EV connected to an EV-PV charger is actively (dis)charging at time $t$.

$$\sum_{v=1}^{v=V} K_{v,c} \leq N_c^{conn} \qquad \forall\, c \qquad (3)$$

$$\sum_{v=1}^{v=V} K_{v,c}\, a_{t,v}^c \leq N_c^{ch} \qquad \forall\, c \qquad (4)$$

Each EV-PV charger feeds $P_{t,c}^{feed}$ or draws $P_{t,c}^{draw}$ power from the EV car park as determined by the EMS. Different EV-PV chargers can exchange power within the car park and these are 'intra-park' power exchanges. When the net 'intra-park' energy exchanges are non-zero, the EV park imports or exports power with the external grid referred to as $P_t^{g(imp)}$, $P_t^{g(exp)}$ respectively.

### B. Trading energy and reserves in the energy market

The ESCo uses the EMS to control the solar powered EV car park for energy trading with the grid. Since $P_t^{g(imp)}$, $P_t^{g(exp)}$ are small relative to the power traded in the market, the ESCo is a price taker and does not influence the market clearing prices. It uses the market settlement point prices $p_t^{e(buy)}$, $p_t^{e(sell)}$ at time instant $t$ for buying and selling power; and reserve capacity prices $p_t^{r(up)}$, $p_t^{r(dn)}$ for up regulation and down regulation for offer of ancillary services. The prices are used as inputs to the EMS for trading (buying and selling) energy between the car park and the grid. Markets like the Electric Reliability Council of Texas (ERCOT) provide different prices for offering capacity reserves for up and down regulation ($p_t^{r(up)} \neq p_t^{r(dn)}$). However, other US markets such as PJM trade up and down regulation as a single product. In order to make the EMS flexible and work with both type of markets, it is designed to take different inputs for $p_t^{r(up)}$ and $p_t^{r(dn)}$ and allow for a requirement that up and down regulation quantities could be equal.

The amount of reserves offered by the EV depends on whether the user enables V2G option or not, i.e. if $x_v^{lb}=0$ or not. When an EV is connected to a bidirectional charger and $x_v^{lb} \neq 0$, even an idle EV that is not charging can offer down regulation and up regulation up to $x_v^{ub}$ and $x_v^{lb}$ respectively.

With a unidirectional charger, an idle EV that is not charging can only offer down regulation up to $x_v^{ub}$.

Power generated by PV panels can be ramped down by moving out of the maximum power point of the PV array. This can be achieved by controlling the DC/DC converter in the EV-PV charger that is connected to the PV array. If $P_{t,c}^{PV}$ is the power generated by the PV at time $t$ at the $c^{th}$ charger, this power can also be offered for down regulation services.

### C. Receding horizon model predictive control

There are two sources of variability in the EV-PV system. The first is the diurnal and seasonal variation in PV generation due to changes in weather. The EMS uses solar forecast information as an input to predict the PV variation. $P_t^{PV(fc)}$ is power generation forecast for an optimally orientated 1kW$_p$ PV array at the car park location with a maximum uncertainty in forecast of $y_{PV}^{fc}$. The second is the variation in the arrival and departure patterns of the EV user and the EV parameters like charging powers limits, efficiency of the battery and SOC.

The EMS is implemented as a receding horizon model predictive control with a time step $\Delta T$ to manage these two variations. The horizon for the model is from 00:00AM to 23:59 PM at midnight. This means for every $\Delta T$ time, the EMS gathers all the inputs in real time, performs the optimization and plans the EV charging for the rest of the day. This means the model inaccuracies with respect to PV forecast or SOC estimation will be corrected by the EMS for every $\Delta T$.

## III. MILP FORMULATION

This section describes the objective function and constraints for the MILP formulation of the EMS. It is important to note that all optimization variables considered are positive.

### A. Acceptance criteria

When an EV arrives at the EV car park, it is connected to one of the *C* number of EV-PV chargers. As mentioned earlier, each EV-PV charger can have up to $N_c^{conn}$ number of EV connected to it. The user links to the EMS and the EMS instructs the user on which EV-PV charger he/she must connect to, based on two 'acceptance criteria'. The first criteria is that the energy demand $d_v$ and parking time, $(T_v^d - T_v^a)$ of all the EVs connected to one EV-PV charger must be such that it is able to meet the energy demand within the stipulated time and within the power limits of the charger, (5). The second criteria is that the arrival SOC of the vehicle must be above the minimum SOC as set by the user (6).

$$\sum_{v=1}^{V}\sum_{c=1}^{C} K_{v,c} \frac{d_v}{T_v^d - T_v^a} \leq Min.\{N_c^{ch} P_c^{EVr}, P_c^{conv}\} \qquad \forall\, v,c \qquad (5)$$

$$B_v^{min} \leq B_v^a \qquad \forall\, v \qquad (6)$$

### B. Constraints: EV and user inputs

The EMS controls the charging power $x_{t,v}^{e+}$ and discharging power $x_{t,v}^{e-}$, up and down regulation reserve capacity $x_{t,v}^{r(up)}, x_{t,v}^{r(dn)}$ of each EV and the power extracted from the PV system $P_{t,c}^{PV}$ of each charger at time *t*. Equations (7) and (8) are used to set the charging power of the EV to zero before the arrival ($t < T_v^a$) and after the departure of the EV ($t \geq T_v^d$).

The binary variable $a_{t,v}^c$ indicates if the EV is connected to the isolated DC/DC converter for charging/discharging and can offer regulation services or not. Since an EV cannot simultaneously charge and discharge, a second binary variable $a_{t,v}^{ch\_v2x}$ is used to ensure that only one of the two variables $x_{t,v}^{e-}$, $x_{t,v}^{e+}$ has a non-zero value for a given $t$. $a_{t,v}^{ch\_v2x}$ is set to 1 for charging and to 0 for V2G. $x_{t,v}^{e-}$, $x_{t,v}^{e+}$ have to be within the power limits of the power converter $P_c^{EVr}$ and the charging and discharging power limits $x_v^{ub}$, $x_v^{lb}$ as set by the EV respectively, as shown in equations (9)-(14).

The maximum charging and discharging powers are also dependent on the SOC of the EV battery as shown in (15) and (16). For example, fast charging of EV battery cannot be done beyond 80% SOC of the battery [28]. Here it is assumed that the maximum charging power linearly reduces from $x_v^{ub}$ to zero when the battery is charged beyond 80% SOC till 100% ($S_{ch}$=0.9). Similarly the maximum discharging power reduces linearly from $x_v^{lb}$ to zero when the battery is discharged below 10% SOC till 0% ($S_{v2x}$=0.1). Even though the dependence of battery power on the SOC is non-linear, this is not considered here as it is beyond the scope of the paper and would prevent us from casting the problem into an MLIP formulantion.

$$x_{t,v}^{e-}, x_{t,v}^{e+}, x_{t,v}^{r(up)}, x_{t,v}^{r(dn)}, a_{t,v}^c = 0 \quad \forall\, t < T_v^a \quad (7)$$
$$x_{t,v}^{e-}, x_{t,v}^{e+}, x_{t,v}^{r(up)}, x_{t,v}^{r(dn)}, a_{t,v}^c = 0 \quad \forall\, t \geq T_v^d \quad (8)$$
$$x_{t,v}^{e+} \leq x_v^{ub}(a_{t,v}^c) \quad \forall\, t,v \quad (9)$$
$$x_{t,v}^{e+} \leq x_v^{ub}(a_{t,v}^{ch-v2x}) \quad \forall\, t,v \quad (10)$$
$$x_{t,v}^{e-} \leq -x_v^{lb}(a_{t,v}^c) \quad \forall\, t,v \quad (11)$$
$$x_{t,v}^{e-} \leq -x_v^{lb}(1 - a_{t,v}^{ch-v2x}) \quad \forall\, t,v \quad (12)$$
$$x_{t,v}^{e-}, x_{t,v}^{e+} \leq P_c^{EVr} \quad \forall\, K_{v,c}=1 \quad (13)$$
$$a_{t,v}^c, a_{t,v}^{ch\_v2x}, a_{t,c}^{d\_f} \in \{0,1\} \quad \forall\, t,c,v \quad (14)$$
$$x_{t,v}^{e+} \leq \frac{-x_v^{ub}}{(1 - S_{ch})}\left(\frac{B_{t,v}}{B_v^{max}} - 1\right) \quad \forall\, t,v \quad (15)$$
$$x_{t,v}^{e-} \leq \frac{-x_v^{lb}}{S_{v2x}}\left(\frac{B_{t,v}}{B_v^{max}}\right) \quad \forall\, t,v \quad (16)$$

Equations (17)-(22) are used to set the initial SOC of the EV battery and estimate the SOC of the battery $B_{t,v}$ based on the charging and discharging efficiency $\eta_v^{ch}$, $\eta_v^{v2x}$ and power $x_{t,v}^{e+}$, $x_{t,v}^{e-}$ respectively. The EMS restricts the SOC to be within the limits $B_v^{min}$, $B_v^{max}$ as set by the EV and/or user. It is assumed that the net energy delivered/absorbed by the EV over one time period due to offer of reserves is zero [12], [13]. Hence, $x_{t,v}^{r(up)}, x_{t,v}^{r(dn)}$ do not appear in (22) for SOC estimation.

$$B_{t,v} = 0 \quad \forall\, t < T_v^a \quad (17)$$
$$B_{t,v} = B_v^a \quad \forall\, t = T_v^a \quad (18)$$
$$B_{t,v} \leq d_v + B_v^a \quad \forall\, t = T_v^d \quad (19)$$
$$B_{t,v} \geq B_v^{min} \quad \forall\, t \geq T_v^a \quad (20)$$
$$B_{t,v} \leq B_v^{max} \quad \forall\, t \geq T_v^a \quad (21)$$
$$B_{t+1,v} = B_{t,v} + \Delta T\left(x_{t,v}^{e+}\eta_v^{ch} - \frac{x_{t,v}^{e-}}{\eta_v^{v2x}}\right) \quad \forall\, t,v \quad (22)$$

*C. Constraints: EV–PV charger and car park*

Under normal operation, the EMS extracts maximum power from the PV array using MPPT as shown in right side of equation (23). The PV power is dependent on the scaling factor $K_c^{PV}$ which scales the installation characteristics (e.g. azimuth, tilt, module parameters) of the PV array connected to the charger $c$ with respect to the 1kWp reference array used for the forecast data $P_t^{PV(fc)}$. The EMS implements PV curtailment if it is uneconomical to draw PV power or if there are distribution network constraints for feeding to the grid. This means that the actual PV power extracted $P_{t,c}^{PV}$ can be lower than the MPPT power of the array, as shown in (23).

The DC-link is used for power exchanges between the three ports of the converter and (24) is the power balance equation for the EV-PV converter. It is assumed that each of the power converters within the EV-PV charger operates with an efficiency $\eta_c^{conv}$. $P_{t,c}^{draw}$, $P_{t,c}^{feed}$ are limited by the power limit of the inverter port $P_c^{conv}$. The binary variable $a_{t,c}^{d\_f}$ is used to ensure that only one of the two variables has a finite value for a given $t$ as shown in (25)-(26).

$$P_{t,c}^{PV} \leq K_c^{PV} P_c^{PVr} P_t^{PV(fc)} \quad \forall\, t,c \quad (23)$$
$$\left\{P_{t,c}^{PV} + P_{t,c}^{draw} + \sum_{v=1}^{v=V}\left(K_{v,c}\, x_{t,v}^{e-}\right)\right\}\eta_c^{conv}$$
$$= \left\{P_{t,c}^{feed} + \sum_{v=1}^{v=V}\left(K_{v,c}\, x_{t,v}^{e+}\right)\right\}/\eta_c^{conv} \quad \forall\, t,c,v \quad (24)$$
$$P_{t,c}^{draw} \leq P_c^{conv}(a_{t,c}^{d\_f}) \quad \forall\, t,c \quad (25)$$
$$P_{t,c}^{feed} \leq P_c^{conv}(1 - a_{t,c}^{d\_f}) \quad \forall\, t,c \quad (26)$$

The intra car-park power exchanges between different EV-PV chargers are related to the power exchanged with the external grid $P_t^{g(imp)}$, $P_t^{g(exp)}$ using (27). They should be within the distribution network capacity $P_t^{DN+}$, $P_t^{DN-}$ as shown in (28)-(29). Both $P_t^{g(imp)}$, $P_t^{g(exp)}$ do not have finite values at the same time because of the way the objective function is formulated and because $p_t^{e(buy)} \geq p_t^{e(sell)}$ at all times

$$\sum_{c=1}^{c=C}\left(P_{t,c}^{draw} - P_{t,c}^{feed}\right) = P_t^{g(imp)} - P_t^{g(exp)} \quad \forall\, t \quad (27)$$
$$P_t^{g(imp)} \leq P_t^{DN+} \quad \forall\, t \quad (28)$$
$$P_t^{g(exp)} \leq P_t^{DN-} \quad \forall\, t \quad (29)$$

Finally, each of the EV offers reserve capacity $x_{t,v}^{r(up)}, x_{t,v}^{r(dn)}$ for up and down regulation. From an EV perspective, the regulation power offered is restricted by the power limitations of the EV ($x_v^{ub}$, $x_v^{lb}$) and the EV charger port $P_c^{EVr}$ as shown in Fig. 3. From the EV-PV charger perspective, the regulation power offered depends on the power rating of the inverter port $P_c^{conv}$ and the power exchanged with the grid $P_{t,c}^{draw}$, $P_{t,c}^{feed}$. This is summarized in equations (30)-(35). While asymmetric reserve offers are assumed here ($x_{t,v}^{r(up)} \neq x_{t,v}^{r(dn)}$), symmetric reserves can be achieved by adding $x_{t,v}^{r(up)} = x_{t,v}^{r(dn)}$ to the constraints.

$$\sum_{v=1}^{v=V} K_{v,c}\, x_{t,v}^{r(up)} + P_{t,c}^{feed} \leq P_c^{conv} \quad \forall\, t,c,v \quad (30)$$

$$\sum_{v=1}^{v=V} K_{v,c}\, x_{t,v}^{r(dn)} + P_{t,c}^{draw} \leq P_c^{conv} \qquad \forall\, t, c, v \quad (31)$$

$$x_{t,v}^{e-} + x_{t,v}^{r(up)} \leq P_c^{EVr}(a_{t,v}^c) \qquad \forall\, K_{v,c}=1 \quad (32)$$

$$x_{t,v}^{e-} + x_{t,v}^{r(up)} \leq x^{ub} \qquad \forall\, t, v \quad (33)$$

$$x_{t,v}^{e+} + x_{t,v}^{r(dn)} \leq P_c^{EVr}(a_{t,v}^c) \qquad \forall\, K_{v,c}=1 \quad (34)$$

$$x_{t,v}^{e+} + x_{t,v}^{r(dn)} \leq -x_v^{lb} \qquad \forall\, t, v \quad (35)$$

### D. Objective function

$$\begin{aligned}
Min.\ C^{opt} =\ & \left(B_v^a + d_v - B_{T_v^d,v}\right) C_v^p \\
& + \Delta T \sum_{t=1}^{T} \left( P_t^{g(imp)} p_t^{e(buy)} - P_t^{g(exp)} p_t^{e(sell)} \right) \\
& - \Delta T\, (1 - y_{PV}^{fc})(\eta_c^{conv})^2 \sum_{t=1}^{T}\sum_{c=1}^{C}\sum_{v=1}^{V} K_{v,c}\{ x_{t,v}^{r(up)} p_t^{r(up)} + x_{t,v}^{r(dn)} p_t^{r(dn)} \} \\
& + \Delta T \sum_{t=1}^{T}\sum_{v=1}^{V} x_{t,v}^{e-}\, C^{V2X} \quad + \quad \Delta T \sum_{t=1}^{T}\sum_{c=1}^{C} P_c^{PVr} P_t^{PV(fc)}\, C^{PV}
\end{aligned} \quad (36)$$

The objective function is to minimize the total costs $C^{opt}$ of EV charging, feeding PV power and offering reserves. The formulation is such that the $C^{opt}$ can be positive or negative. It has five components respectively:

- The penalty to be paid to the user if the energy demand $d_v$ is not met by the departure time $T_v^d$. $C_v^p$ is EV user specific and the penalty can be different for each user based on EV battery size, tariff policy and customer 'loyalty' program.
- The cost of buying and selling energy from the grid based on the settlement point prices $p_t^{e(buy)}, p_t^{e(sell)}$. The market dynamics will ensure that $p_t^{e(buy)} \geq p_t^{e(sell)}$
- Income $S^{as}$ obtained from offering reserve capacity $x_{t,v}^{r(up)}, x_{t,v}^{r(dn)}$ to the ISO. $(\eta_c^{conv})^2$ indicates the energy losses in the two step conversion between the EV and grid port of the EV-PV charger. Since the reserves offered to the grid have to be guaranteed and the uncertainty in the PV forecast is $y_{PV}^{fc}$, only a fraction $(1 - y_{PV}^{fc})$ of the available reserves are guaranteed and sold to the ISO.
- EV battery capacity degrades due to the additional life cycles caused by the V2G operation and EV user is compensated for this loss. Typical value of $C^{V2X}$=4.2¢/kWh based on analysis in [29], [30].
- PV power that is used to charge the EV need not always be free of cost. If the PV is installed by a third-party, it can be obtained at a pre-determined contractual cost of $C^{PV}$.

### E. MILP implementation

The EMS engine is implemented in C# leveraging Microsoft Solver Foundation for algebraic modeling in Optimization Modeling Language (OML). MS SQL Server database is used to warehouse system inputs, namely the EV, charger, network and market data as well as the decision outputs that is sent to the EV-PVs in the field. The MILP formulation is solved using branch-and-bound (B&B) algorithm using 'LPsolve' open source solver. One of the main advantages of the B&B algorithm is that, given enough computation time, it guarantees global optimality despite the non-convex nature of the problem. The EV-PV chargers will be interfaced with the output database to implement the optimal power profiles.

## IV. SIMULATION RESULTS

Simulations are performed to test the validity of the proposed MILP formulation and to quantify the reduction in costs of EV charging from PV with respect to AR and IMM.

### A. Simulation parameters

Settlement point prices (SPP) and prices for reserve capacity (REGUP, REGDN) are obtained from the ERCOT day-ahead market (DAM) for Austin, Texas for 2014 for load zone LZ_AEN. These are wholesale energy prices with a data resolution of 1hr. Since separate values for $p_t^{e(sell)}$ was not available, it is assumed that $p_t^{e(sell)}$=0.98*$p_t^{e(buy)}$.

For 2014, the largest values observed for $p_t^{e(buy)}$, $p_t^{r(up)}$, $p_t^{r(dn)}$ were 136.47¢/kWh, 499.9¢/kWh and 31¢/kWh respectively while the average values were 3.9¢/kWh, 1.25¢/kWh, 0.973¢/kWh. It can be clearly seen than energy prices are normally much higher than regulation prices, but there are several instances where it is otherwise.

The PV generation data is obtained from the Pecan Street Project database for a house in the Mueller neighborhood with an 11.1 kW PV system [31]. The data resolution is 1min. The power output is scaled down for a 1kW system for use as $P_t^{PV(fc)}$ with $y^{PV(fc)}$=10%. It is assumed that the PV installation at the car park is owned by the workplace and hence $C^{PV}$=0

The EV arrival and departure times and SOC requirements are listed in TABLE I for 6 EVs. The EV data imitates the capacity of a Tesla Model S, BMW i3 and a Nissan Leaf. For all the EVs, $B_v^{min}$=5kWh, $x_v^{ub}$=50kW, $x_v^{lb}$=(-10kW), $\eta_v^{ch} = \eta_v^{v2x}$=0.95, $C_v^p$=1$/kWh, $C^{V2X}$= 4.2¢/kWh. The penalty $C_v^p$ is approximately 25 times the average wholesale ERCOT electricity price of 3.9¢/kWh.

There are 4 EV-PV chargers and the TABLE I shows the connections of the 6 EVs to the 4 chargers in '*Chr conn.*'. 10kW$_p$ PV is connected to each of chargers 1,2,4 and no PV is connected to charger 3. Chargers 1,4 have two EV connected to them. $N_c^{ch}$=1 for all chargers, which means that only one of the two EV can be charged at a time for chargers 1,4. The following parameters are used: $\eta_c^{conv}$=0.96, $P_c^{EVr}$=$P_c^{conv}$=10 kW. $\Delta T$=15min for all simulation. $P_t^{DN+} = P_t^{DN+}$ =40kW.

### B. Simulation results

#### 1) Average rate and immediate charging

The net costs of EV charging and PV sales for average rate

TABLE I
EV AND EV-PV CHARGER DATA

| v | $T_v^a$ | $T_v^d$ | d | $B_v^a$ | $B_v^{max}$ | Chr conn. | $P_c^{EVr}$ $P_c^{conv}$ |
|---|---|---|---|---|---|---|---|
| | (h) | | | (kWh) | | | (kW) |
| 1 | 900 | 1700 | 40 | 20 | 85 | 1 | 10 |
| 2 | 830 | 1630 | 30 | 20 | 60 | 1 | 10 |
| 3 | 930 | 1730 | 10 | 5 | 24 | 2 | 10 |
| 4 | 900 | 1700 | 40 | 20 | 85 | 3 | 10 |
| 5 | 830 | 1630 | 30 | 20 | 60 | 4 | 10 |
| 6 | 930 | 1730 | 10 | 5 | 24 | 4 | 10 |

$C^{ar}$ and immediate charging $C^{imm}$ is estimated using (1), (37).

$$C^{ar}, C^{imm} = C^{ev} - S^{PV}$$
$$= \Delta T \sum_{t=1}^{T} \sum_{v=1}^{v=V} (\eta_c^{conv})^2 x_{t,v}^{e+} p_t^{e(buy)} -$$
$$\Delta T \sum_{t=1}^{T} \sum_{c=1}^{C} (\eta_c^{conv})^2 P_c^{PVr} P_c^{PV(fc)} \left( p_t^{e(sell)} - C^{PV} \right) \quad (37)$$

where $C^{ev}$ is the EV charging costs and $S^{PV}$ the revenues from PV sales. For AR, $x_{t,v}^{e+} = x_v^{e(ar)}$ and for IMM, $x_{t,v}^{e+} = P_c^{EVr}$. With AR and IMM, there is no provision to provide V2G, regulation services or multiplexing of chargers due to the absence of communication with an EMS. The peak power for the car park for IMM would be 60kW and 20kW for AR charging for 6EV based on (1).

Fig. 4 and TABLE II shows the net costs $C^{ar}, C^{imm}$ estimated for 2014 with the corresponding mean and standard deviation (SD). Three vital observations can be made. First, there is a large variation in net costs, ranging between {1.35$, 24.17$} and {-19.58$, 40.43$} for AR and IMM respectively. This is mainly due to the varying energy prices in ERCOT. The costs went negative for IMM on certain days indicating that the ESCo got paid by the ISO! It must be remembered that PV sales $S^{PV}$ for both strategies is the same as shown in TABLE II. Second, IMM charging was found to be better than AR in summer and vice versa in winter, with IMM charging net costs being cheaper than AR for 233 days. Third, the average net cost per day for 2014 for AR and IMM was found to be 3.79$ and 2.90$, with IMM being cheaper than AR by 31.7%. This is because, EVs are charged in morning for IMM when ERCOT prices are generally lower when compared to prices in the afternoon.

*2) Optimized charging costs*

Using the MILP formulation of the optimized charging (OPT) described in section III, the net costs $C^{opt}$ are determined for each day of 2014 and shown in Fig. 4 and TABLE II. The benefits of the MILP optimization can be clearly seen in the figure, where the optimized net costs are much lower than IMM and AR. $C^{opt}$ range is {-42.91$, 11.56$}, which is much lower than IMM and AR. Due to the large penalty $C_v^p$=1$/kWh, EVs were always charged up to the required departure SOC.

EV charging costs $C^{ev}$ (not net cost!) are estimated separately for AR, IMM and OPT and shown in TABLE II. It can be seen that mean value of $C^{ev}$ is not that different between IMM and OPT. The reason is that the objective function is not optimized to reduce charging costs alone but increase the sale of PV power and reserves as well.

The percentage reduction in costs $C_\%^{imm}$, $C_\%^{EV-PV}$ is estimated based on AR charging costs $C^{ar}$ using (38)-(39) and shown in Fig. 5. $C^{ar}$ was chosen as a reference as the costs never go negative and don't have values close to zero.

$$C_\%^{imm} = 100(C^{ar} - C^{imm})/C^{ar} \quad (38)$$
$$C_\%^{opt} = 100(C^{ar} - C^{opt})/C^{ar} \quad (39)$$

As can be seen, the proposed optimized charging results in a cost reduction $C_\%^{opt}$ in the range of 31.74% to 650.81%, with a mean of 158.63% with respect to AR charging. A reduction of >100% results in the net cost to be negative. This means that the EV car park receives money for the EV charging and sale of PV and reserves rather than having to pay overall.

The large cost reduction is a result of aggregating the multi-aspect PV and EV problems into a single MILP formulation. This results in the sale of PV and V2G power when prices are high, buying of EV charging power when prices are low and continuous sale of ancillary services. The current MILP formulation is such that IMM or AR will be a special case of optimized charging as dictated by the PV forecast and market prices. Further, the sharing of a single charger to charge several EVs results in a reduction of charging infrastructure cost. While these costs have not be included in the estimate, they can be up to 15,000$ for 10kW chargers with $N_c^{conn}$=4.

MILP solve times were in the range of 11.2-17.3s with a relative MILP gap of 0.015%. The mean solve-time was 13.05s with a standard deviation of 1.09s. A Windows PC with Intel Xeon 2.4Ghz CPU and 12GB RAM was employed.

It must be kept in mind that even though wholesale DAM prices and small EV fleet have been used in this simulation, the formulation is generic to be used with large EV fleet, real-time market (RTM) and retail electricity prices as well.

## V. CONCLUSIONS

EV charging from PV can be controlled to achieve several motives – to take advantage of time of use tariffs, provide ancillary services or follow the PV production. However, the common approach is that each of these applications are solved as separate optimization problems which leads to several EV charging profiles. This is impractical, as a single EV cannot be controlled at the same time with different charging profiles.

TABLE II
CHARGING COSTS, PV SALES AND NET COSTS - MEAN, SD ($)

| [Mean, SD] | AR | IMM | OPT |
|---|---|---|---|
| $S^{PV}$ | 4.41, 2.81 | 4.41, 2.81 | - |
| $C^{ev}$ | 8.21, 3.21 | 7.32, 3.87 | 7.30, 1.92 |
| $C^{ar}, C^{imm}, C^{opt}$ | 3.79, 2.13 | 2.90, 42.01 | -1.53, 3.92 |
| $C_\%^{imm}, C_\%^{opt}$ (%) |  | 31.72, 61.26 | 158.63, 87.88 |

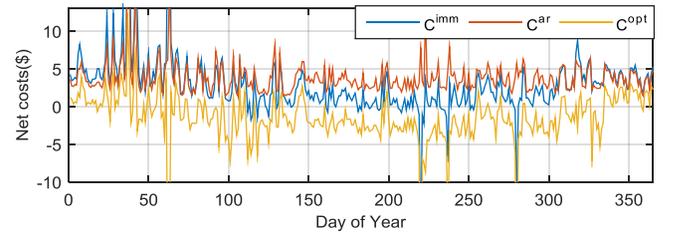

Fig. 4. Cost of charging the EV fleet by average rate, immediate and the proposed optimized charging strategy (top); zoomed view (bottom)

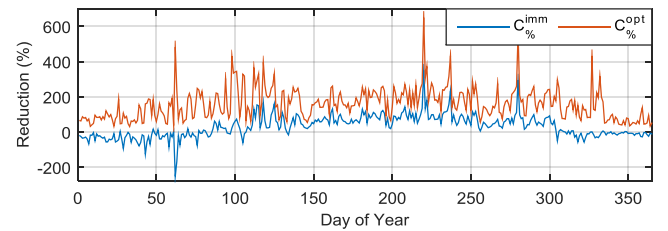

Fig. 5. Percentage reduction in the charging cost for the proposed charging strategy and immediate charging with respect to average rate charging.

Hence it is vital to make a single problem formulation that bundles several applications together so that one optimal EV charging profile is obtained.

In this paper, an MILP formulation has been proposed for charging of an EV fleet from PV that has several application built into one - charging of EV from PV, using time of use tariffs to sell PV power and charge EV from the grid, implementation of V2G for grid support, using EV to offer ancillary services in the form of reserves and considering distribution network capacity constraints. The scheduling of the connection of a single EVSE to several EV has been formulated for the first time in this work. This provides the ability to use lower capacity EVSE at workplaces resulting in substantial reduction in the cost of EV infrastructure.

The MILP optimization has been implemented as a receding horizon model predictive control and operates with a fixed time period. Using 2014 data from Pecan Street Project and ERCOT market, simulations were performed for an EV fleet of six connected to four chargers. The formulation of five applications into one resulted in large savings in the range of 31.74% to 650.81% with respect to average rate charging. The MILP formulation is generic and can be adapted to different energy and ancillary markets, EV types, PV array installations and EVSE.